\documentclass[11pt]{amsart}
\usepackage{graphicx}
\usepackage{amssymb}

\newcommand{\R}{\mathbb{R}}
\newcommand{\half}{\frac{1}{2}}

\newcommand{\quarter}{\frac{1}{4}}
\newcommand{\twoovers}{\frac{2}{s}}
\newcommand{\sQ}{\frac{1}{5} + \frac{1}{5} \sqrt{11}}

\def\C{{\mathbb C}} 
\def\R{{{\mathbb R}}}

\newcommand{\tstar}{T^*}

\newcommand{\al}{\alpha}

\newcommand{\ep}{\epsilon}

\newcommand{\proj}{\mathbb{P}}

\newtheorem{theorem}{Theorem}
\newtheorem{remark}[theorem]{Remark}

\newtheorem{proposition}[theorem]{Proposition}
\newtheorem{lemma}[theorem]{Lemma}
\newtheorem{corollary}[theorem]{Corollary}

\numberwithin{equation}{section}
\numberwithin{theorem}{section}

\title[Ground state mass concentration
for NLS below $H^1$]{Ground state mass concentration in the $L^2$-critical nonlinear
Schr\"odinger equation below $H^1$}
\author{J.~Colliander}
\address{University of Toronto and The Fields Institute}
\thanks{J.C. is partially  supported  by N.S.E.R.C. Grant RGPIN 250233-03 and the Sloan Foundation.}
\author{S.~Raynor}
\address{University of Toronto and The Fields Institute}
\curraddr{Wake Forest University}
\author{C.~Sulem}
\address{University of Toronto and The Fields Institute}
 \thanks{C.S. is partially supported by 
            N.S.E.R.C. Grant  OGP 0046179.}
\author{J.~D.~Wright}
\address{McMaster University and The Fields Institute}
\curraddr{University of Minnesota}


\begin{document}
\date{\today}

\begin{abstract}
We consider finite time blowup solutions of the
$L^2$-critical cubic focusing
 nonlinear Schr\"odinger  equation on $\R^2$. 
 Such functions, when in $H^1$, are known  to concentrate a fixed
 $L^2$-mass (the mass of the ground state) at the point of blowup.
 Blowup solutions from initial data that is only in $L^2$ are
 known to concentrate at least a small amount of mass. 
In this paper we consider the intermediate case of
 blowup solutions from initial data in $H^s$, with $1 > s > s_Q$, 
where $s_Q \le \sQ$.  Our main
 result is that such solutions, when radially symmetric,
 concentrate at least the mass of
 the ground state at the origin at blowup time.

\end{abstract}

\maketitle


\section{Introduction}

Special interest has recently been devoted to the existence and long-time
behavior of solutions with low regularity 
to nonlinear Schr\"odinger equations. These questions 
were mainly investigated for defocusing\footnote{Low regularity
  global well-posedness results have also been obtained for other
  Hamiltonian evolution equations. Orbital instability
  properties of solitons subject to rough perturbations for focusing NLS equations have been
  studied as well \cite{CKSTT:StabI}, \cite{CKSTT:StabII}. }    
equations with a global-in-time a priori $H^1$ upper bound \cite{BRefine}
 \cite{CKSTT:MRL} \cite{Tzir} \cite{CKSTT:CPAM}. 
In this
article, we are interested in a detailed description of rough solutions
in a situation where the $H^1$-norm of certain smooth solutions 
blows up in a finite time.

We consider the initial value problem for the two-dimensional,
cubic, focusing nonlinear Schr\"odinger (NLS) equation:
\begin{equation}
\label{nls}
\left\{
\begin{matrix}
i u_t + \Delta u  + |u|^2 u = 0,  \\
 u(0,x) = u_0(x),& x \in \R^2,
\end{matrix}
\right.
\end{equation}
which is $L^2$-critical. 
This refers to the property that both the equation and the $L^2$-norm 
of the solution are invariant under the scaling transformation
$u(t,x) \to u_\lambda(t,x) = \lambda u(\lambda^2 t,  \lambda x)$.
This problem is locally well-posed\footnote{For
  $s=0$, the size of the interval of existence depends upon the initial
  profile $u_0$ \cite{CWzero};
for $s > 0$, the $H^s$ norm of the data determines
  the size of the existence interval.}  for initial data 
in $H^s$ with $s\ge0$ \cite{CW}.
We recall that the following quantities, if finite, are conserved:
\begin{equation*}
\begin{split}
Mass&=M[u(t)] = \|u(t)\|_{L^2}^2,\\
Energy&=E[u(t)] = \half \|\nabla u(t)\|_{L^2}^2 - \quarter \|u(t)\|_{L^4}^4.
\end{split}
\end{equation*}
We will frequently refer to $\half \|\nabla u(t)\|_{L^2}^2$ as the {\it kinetic energy} of
the solution, and $E[u(t)]$ as the {\it total energy}.

It is known that
there exist explicit finite time blowup
solutions to \eqref{nls}; for sufficiently smooth and decaying
initial data, the virial identity provides a sufficient condition
guaranteeing that finite time blowup occurs.  For a solution which
blows up in finite time, let $[0, T^*)$ be the maximal (forward) time interval
  of existence.

A specific property of critical collapse is the phenomenon of mass
concentration, often referred to in the physical literature as
{\it{strong collapse}} \cite{ZK} (see also \cite{SuSu} for a review): 
 $H^1$-solutions concentrate a finite amount of mass 
 in a neighborhood of the focus of 
width slightly larger than $(T^*-t)^{1/2}$. 
Heuristic arguments suggest
that this phenomenon does not occur in supercritical
nonlinear Schr\"odinger blowup. For $H^1$-solutions of \eqref{nls},
there is a precise lower bound on the amount of concentrated mass,
namely the mass of the ground state
$Q$ \cite{T}, \cite{MT1}, where $Q$ is the
 unique positive solution (up to translations) of
\begin{equation}
  \label{groundstateeq}
  \Delta w -w + |w|^2 w = 0.
\end{equation}
In addition to the scaling properties of the NLS equation, 
the main ingredients in the proof that $H^1$ blowup solutions
concentrate at least the mass of the ground state are:
(i) the conservation of the energy, (ii) a precise Gagliardo-Nirenberg
inequality \cite{W0} which implies that
nonzero $H^1$-functions of non-positive energy
have at least ground state mass.

The purpose of this work is to address the phenomenon of mass
concentration in the spaces $H^s$, $s<1$, where the conservation
of energy cannot be used. In the setting of merely
    $L^2$ initial data, if global well-posedness fails to hold 
for \eqref{nls}
({\it i.e.} $T^* < \infty$), then a
nontrivial parabolic concentration of $L^2$-mass occurs \cite{BRefine}
 as $t \uparrow T^*$:
\begin{equation}
  \label{tinyconc}
  \limsup_{t \uparrow T^* }  \sup\limits_{\begin{matrix} {\mbox{cubes}}~ I \subset \R^2 \\
                                           {\mbox{side} (I) < (T^* -
                                           t)^\half} \end{matrix}}
\left( \int_I |u(t,x)|^2 dx \right)^\half \geq \eta ( \|u_0
\|_{L^2} )
> 0.
\end{equation}
Unlike the $H^1$-case, there is no explicit quantification
 on the lower bound $\eta ( \|u_0 \|_{L^2})$.

A natural question{\footnote{{S. Keraani has announced that the lower
      bound in \eqref{tinyconc} may be taken to be a fixed constant $\delta_0$
      independent of the initial size in $L^2$. }}}, highlighted in
\cite{MV}, is to determine
whether tiny $L^2$-mass concentrations can occur when $u_0 \in
L^2$. The conjectured answer is no. Solutions of \eqref{nls} with
a finite maximal (forward) existence interval are expected to
concentrate at least the $L^2$-mass of the ground state.
Our main result corroborates this expectation,
at least for $H^s$-solutions with $s$ just below 1.
\begin{theorem}
\label{main result} There exists $s_Q \le \sQ$ such that
the following is true for any $s > s_Q $. Suppose $H^s \ni u_0
\longmapsto u(t)$ solves \eqref{nls} on the maximal (forward) time
interval $[0,T^*)$, with $T^*  < \infty$. Moreover, assume that $u_0$
is radially symmetric.\footnote{The
  non-radial case is amenable to treatment by
 employing the methods of  compensated compactness as used in
 \cite{W},  \cite{Nawa} (see also \cite{CazenaveCourant}).}
Then for any positive $\gamma (z) \uparrow \infty$ arbitrarily slowly as
$z \downarrow 0$ we have
\begin{equation}
\label{sparabolic}
\limsup_{t \uparrow \tstar} \|u\|_{L^2_{\{|x| < (T^* - t)^{s/2}\gamma (T^* - t) ~\}}} \ge \|Q\|_{L^2}.
\end{equation}
\end{theorem}

The proof consists of an imitation of the $H^1$ argument (as
presented in \cite{CazenaveCourant}) with the energy, which is
infinite in the $H^s$-setting, replaced by a modified energy
introduced in \cite{CKSTT:MRL}. The idea is to apply to 
the $H^s$-solution a smoothing operator to make it  $H^1$ and define
the usual energy of this new object.  The crucial point here is to prove
that the modified total energy grows more slowly than the modified
kinetic energy.
In Section 2, we  prove Theorem \ref{main result} assuming
Proposition \ref{umbrella} and Corollary
\ref{modifiedcwlowercor} stated below. Proposition \ref{umbrella} contains
a key  upper bound  for the modified energy in terms of the
$H^s$-norm of the solution.  Its  proof,
 which relies upon the local-in-time
theory developed in Section \ref{localtheorysection},
 is postponed to   Section
\ref{umbrella lemma section}. 

\begin{remark}
As mentioned above, the expectation is that parabolic concentration of at least
the ground
state mass holds true even for $L^2$ initial data.
However, the analysis showing that the
modified total energy grows more slowly than the modified kinetic
energy (quantified by the statement that $p(s) < 2$ in Proposition
\ref{umbrella}), requires taking $s$ close to 1. Also, note that the
concentration width $(T^* -t )^{s/2 +}$ obtained in \eqref{sparabolic}
is much larger than $(T^* -t)^{1/2 +}$ 
within which ground
state mass concentration is expected to occur.
\end{remark}

\begin{remark}
The two-dimensional nature of the analysis only appears directly in the use of
the Sobolev inequality in the proof of Lemma \ref{6linearlemma}.  In
higher dimensions the power of the nonlinearity of the $L^2$-critical
problem is non-integer, which would make the multilinear analysis of
the modified energy more technical. 
\end{remark}

\begin{remark}
Our methods only enable us to control the modified energy at time $t$
by the supremum up to time $t$ of the modified kinetic energy.  This
is the source of the $\limsup$ in Theorem \ref{main result}, which
does not appear in the $H^1$ result in \cite{T}.
  Because there is
no known upper bound for the rate of blowup of the $H^s$ norm, we are unable
to control the relationship between times  at which the blowup is occurring
maximally and those where it is occurring minimally.  Any such result
would, in addition
to its inherent interest, allow us to strengthen our result to one 
with a $\liminf$.
For further discussion of the monotonicity of the blowup, but in the
context of the critical generalized KdV equation, see pp. 621--623
of  \cite{MartelMerle:JAMS}.
\end{remark}

\begin{remark}
We use $C$ to denote constants which do not depend on the 
crucial parameters and
variables.  Such constants will {\it not} depend on the high-frequency
cut-off parameter (denoted $N$) or the time $t$, but may
depend on $s$, $\tstar$, or  $\|u_0\|_{H^s}$.
\end{remark}

In the proof of Theorem \ref{main result},
nowhere do we make the hypothesis that $\|u_0\|_{L^2}\ge\|Q\|_{L^2}$.
Thus, since solutions conserve mass and the concentration is shown
to be in excess of the ground state mass, we prove the following 
corollary about the global well-posedness of \eqref{nls}.
\begin{corollary}
\label{gwpcor}
There exists $s_Q \leq \frac{1}{5} + \frac{1}{5} \sqrt{11}$ such
that, if $u_0 \in H^s, ~s > s_Q$ is radially symmetric 
and $\| u_0 \|_{L^2} < \| Q \|_{L^2}$ then the initial value problem 
\eqref{nls} is globally well-posed.
\end{corollary}

\begin{remark}\footnote{This remark emerged in correspondence among
the authors of \cite{CKSTT:MRL}.}
A statement stronger than Corollary \ref{gwpcor} may be
inferred from earlier work.
The analysis in \cite{CKSTT:MRL} of the defocusing analog of
\eqref{nls} relies upon the local well-posedness theory and two other 
inputs: the almost conservation of the modified energy (Proposition
3.1 of \cite{CKSTT:MRL}) and the obvious fact that, in the defocusing
    case, the modified total
energy controls the modified kinetic energy. The almost conservation
property depends upon the local-in-time space-time boundedness properties
of the solution and also holds in the focusing case. Since
the smoothing operator $I_N$ (see \eqref{INdefined} below)  shrinks the
$L^2$ size of functions, it can be shown that the modified total
energy does indeed control the modified kinetic energy provided that
the initial data satisfies the smallness condition $\| u_0 \|_{L^2} <
\|Q \|_{L^2}$. Thus, Corollary \ref{gwpcor} actually holds without the
radial symmetry restriction and for $s> \frac{4}{7}$. 
\end{remark}

{\emph{Acknowledgments:}} We thank P. Raphael, D. Slepcev
and W. Staubach for interesting
conversations related to this work. We gratefully acknowledge support
from the Fields Institute where this research was carried out.

\section{Mass Concentration}

We wish to find a replacement for energy conservation.  To do so, we
define 
the smoothing operators $I_N :H^s \to H^1$ used in \cite{CKSTT:MRL}:
\begin{equation} 
\label{INdefined}
\widehat{I_N u}(\xi) = m(\xi) \hat{u}(\xi)
\end{equation}
where
\begin{equation}
\label{mdefined}
m(\xi) = \left\{ \begin{matrix}
1,  & |\xi | \leq N \\
(\frac{|\xi|}{N})^{s-1},  & |\xi | > 3 N
\end{matrix}
\right.
\end{equation}
with $m(\xi)$ smooth, radial,  and monotone in $|\xi|$.  Note that we will sometimes drop the
subscript $N$ when that will not lead to confusion. The following properties of $I_N$
are easily verified:
\begin{equation}
\label{easy I properties}
\begin{split}
&\|I_N u\|_{L^2}  \le  \|u\|_{L^2}\\
\|u\|_{H^s} \le &\|I_N \langle D \rangle u\|_{L^2}  \le  N^{1-s} \|u\|_{H^s},
\end{split}
\end{equation}
where $\langle D \rangle$ is the multiplier operator with symbol $(1 + |\xi|^2)^{1/2}$.
We define the {\it blowup parameters}
associated to the $H^s$-norm of the solution
\begin{equation}
\label{lambdadefined}
\lambda(t)=\|u(t)\|_{H^s}
\end{equation}
and
\begin{equation}
\label{Lambdadefined}
\Lambda(t) = \sup_{0\le \tau \le t} \lambda (\tau) .
\end{equation}
We also define the {\it{modified blowup parameters}}
\begin{equation}
  \label{sigmadefined}
  \sigma(t) = \| I_N \langle D \rangle u(t) \|_{L^2}
\end{equation}
and
\begin{equation}
  \label{Sigmadefined}
  \Sigma(t) = \sup_{0 \leq \tau \leq t} \sigma(\tau).
\end{equation}
For finite time blowup solutions $u$, 
the blowup parameters $\sigma(t)$ and $\lambda(t)$ tend to $\infty$ as $t
    \rightarrow T^*$, and can be compared using
\eqref{easy I properties}.

We will exploit the freedom to choose $N$ very large. Indeed, for
fixed $T < T^*$, we will show that $N= N(T)$ may be chosen so large
that the modified energy
 $E[I_N u(T)]$ is much smaller than the modified kinetic
energy $\| I_N \nabla u (T) \|_{L^2}^2.$

\begin{proposition}
  \label{umbrella}There exists
$s_Q \leq \frac{1}{5} + \frac{1}{5} \sqrt{11}$ such that for all $s> s_Q$ there exists
$p(s) < 2$ with the following holding true: If $H^s \ni u_0
\longmapsto u(t)$ solves \eqref{nls} on a maximal (forward) finite
existence interval $[0,T^*)$ then for all $T < T^*$ there exists
$N=N(T)$ such that
\begin{equation}
  \label{energycontrol}
  |E[I_{N(T)} u(T)]| \leq C_0 (\Lambda(T))^{p(s)}
\end{equation}
with $C_0 = C_0 (s, T^*,  \| u_0 \|_{H^s} )$. Moreover, $N(T) = C
(\Lambda(T))^\frac{p(s)}{2(1-s)}$.
\end{proposition}

\begin{remark}  For $s\geq 1$ we set $p=0$.
  In this case, Proposition \ref{umbrella} is reduced to the statement
  that the energy remains bounded, which is true since 
  it is conserved.
  Note also that we have chosen $N =
  N(\Lambda)$ so the time dependence of $N$ only enters through that of
the blowup parameter
   $\Lambda =\Lambda (t).$
\end{remark}

Proposition \ref{umbrella} gives a control on the growth of the modified
energy as $t$ approaches $T^*$. During the finite time blowup
evolution, the modified kinetic energy explodes to
infinity. The above Proposition shows that the total 
modified energy grows more slowly than its kinetic component.
It is the key element in the proof of Theorem \ref{main result}, requiring
somewhat delicate harmonic analysis  estimates to prove an almost
conservation property of the modified energy $E[I_N u ]$ in terms of
space-time control given by local existence theory. For the sake of a 
clear presentation,  we start by  proving  Theorem \ref{main result}
 assuming Proposition \ref{umbrella} and a lower bound on the rate of
 blowup of $\| I_N \langle D \rangle u \|_{L^2}$ expressed in Corollary
 \ref{modifiedcwlowercor}. 

\begin{proof}[Proof of Theorem \ref{main result}] 
It is carried out  in four steps.

{\bf{a. Rescaling and weak convergence.}} \newline
Let $\{t_n\}_{n=1}^\infty$ be a
sequence such that $t_n \uparrow \tstar$
and for each $t_n$
$$
\|u(t_n)\|_{H^s}=\Lambda(t_n).
$$
We call such a sequence  {\it maximizing} and denote $u_n = u(t_n)$.
Define
$$
I_N u_n = I_{N(t_n)} u(t_n)
$$
with $N(t_n)$ taken{\footnote{The reason for this choice is made clear
    at equation \eqref{Nchoice}.}} as in the Proposition  \ref{umbrella}.  We rescale these
as follows:
\begin{equation}
\label{rescale}
v_n(y) = \frac{1}{\sigma_n} I_N u_n(\frac{y}{\sigma_n})
\end{equation}
where
$$
\sigma_n = \| I_N \langle D \rangle u_n \|_{L^2} = \sigma (t_n ).
$$
We have, from \eqref{easy I properties}, conservation of the $L^2$
norm 
and $\Lambda_n \rightarrow \infty$, that
\begin{equation}
\label{lower bound}
\Lambda(t_n) \le \sigma_n
\end{equation}
along a maximizing sequence so $\sigma_n \to \infty$ as $n\to\infty$.
The lower bound of $\sigma_n$ by $\Lambda(t_n)$ does not necessarily
hold for arbitrary sequences, but does so along maximizing
sequences.

The rescaling \eqref{rescale} leaves $L^2$-norms unchanged.
Since $\|I_N \cdot\|_{L^2} \le \| \cdot\|_{L^2}$ for all $N$
and since mass is conserved, we have
$\|v_n\|_{L^2}\le\|u_0\|_{L^2}$ uniformly in $n$.
By our choice of $\sigma_n$,
we have $\|\nabla v_n\|_{L^2} \le 1$ for all $n$.  In fact, by
construction,  we have
\begin{equation}
\label{non zero}
\lim_{n\to\infty} \|\nabla v_n\|_{L^2} =1 .
\end{equation}
Thus, $\{v_n\}$ is a bounded sequence in $H^1$ and has a weakly
convergent subsequence, which we also call $\{v_n\}$.  There
exists an asymptotic object $v \in H^1$ such that
$$v_n \rightharpoonup v$$ in $H^1$.  

{\bf{b. Compactness and energy of the rescaled asymptotic object}}.

Since $u$ is assumed to be radially symmetric, so are the
$v_n$, and we can apply the following Lemma from \cite{Strauss}.
\begin{lemma}
({\bf Radial compactness lemma})  If $\{f_n\} \subset H^1 (\R^2) $ is a bounded
sequence of radially symmetric functions, then there exists a subsequence 
(also denoted $\{f_n\}$) and a function $f \in H^1 (\R^2) $ such that for all
 $2 < q < \infty$, $$f_n \to f$$ in $L^q$.
\end{lemma}

Thus,
$$
v_n \to v
$$
in $L^4$. This strong $L^4$-convergence is important due to the
appearance of the $L^4$ norm in the energy, which we now examine.  
The only usage of the radial symmetry assumption in the argument is to
obtain compactness of $\{ v_n \}$ in $L^4$.

We have, from Proposition \ref{umbrella} and (\ref{lower bound}),
\begin{equation*}
\begin{split}
|E[v_n] | = & \sigma_n^{-2} |E[I_N u_n]|\\
      \le& C \sigma_n^{-2} \Lambda^{p(s)}(t_n)\\
       \le & C \Lambda^{p(s)-2}(t_n).
\end{split}
\end{equation*}
Since $s > s_Q $, we have $p(s) < 2$, so
$$
|E[v_n] |\to 0
$$
as $n \to \infty$.

The fact that  the energy of the functions  $v_n$ goes to zero is
 useful in two important
ways.  First, by (\ref{non zero}) and the strong $L^4$-convergence,
\begin{equation}
\begin{split}
0 & = \lim_{n\to\infty} |E[v_n] |\\
  & = \lim_{n\to\infty} \left|  \half \|\nabla v_n\|_{L^2}^2 -
  \quarter \|v_n\|_{L^4}^4  \right|\\
  & = \left|\half - \quarter \|v\|_{L^4}^4 \right|.
\end{split}
\end{equation}
This ensures that $v \ne 0$. Second, the $L^2$ norm is a lower
semi-continuous function for weakly convergent sequences, so
$$
0=\lim_{n\to\infty}E[v_n] \ge E[v].
$$

{\bf{c. Non-positive energy implies at least ground state mass.}}

The fact that nonzero functions of nonpositive energy have at least the mass of the ground state is a consequence of the Gagliardo--Nirenberg inequality:
\begin{equation}
\| w\|_{L^4}^4 \ \leq 
 \ C_{opt}
\| w \|_{L^2}^{2}  \|\nabla w\|_{L^2}^2.
\label{Maj}
\end{equation}
The optimal constant, obtained by minimizing the functional
\begin{equation}
J(f) = \frac{\|\nabla f\|_{L^2}^{2}
 \|f\|_{L^2}^2}{ \|f\|_{L^4}^4 }
\label{JJeq}
\end{equation}
among all functions  $f  \in H^1(\R^2)$, 
is found to be $C_{opt} =2/{\|Q\|_{L^2}^2}$ \cite{W0}.

Thus, the asymptotic object $v$ satisfies
\begin{equation}
\label{asymptotic object has gsm}
\|v\|_{L^2} \ge \|Q\|_{L^2}.
\end{equation}

{\bf{d. Scaling back to the original variables.}}

Now we complete the proof.  To prove \eqref{sparabolic}, it suffices to show for any
$\ep>0$ that
$$
\lim_{n\to\infty} \|u(t_n)\|_{L^2_{\{|x|<(T^* - t_n )^{s/2} \gamma(T^*
    - t_n)\}}}  > \|Q\|_{L^2} -
\ep .$$
  Fix $\ep>0$.  Since $N(t_n)$ goes to
$\infty$, we have
$$
\lim_{n\to\infty} \|u(t_n)\|_{L^2_{\{|x|<(T^* - t_n )^{s/2} \gamma(T^*
    - t_n)\}}} =
\lim_{n\to\infty}\|I_{N(t_n)} u(t_n)\|_{L^2_{\{|x|<(T^* - t_n )^{s/2} \gamma(T^*
    - t_n)\}}}.
$$
Recalling the definition of the functions $v_n$, we have for all $n$
$$
\|I_{N(t_n )}u(t_n)\|_{L^2_{\{|x|<(T^* - t_n )^{s/2} \gamma(T^*
    - t_n) \}}} = \|v_n\|_{L^2_{\{|y|<(T^* - t_n )^{s/2} \gamma(T^*
    - t_n)
\sigma_n \}}}.
$$
By (\ref{asymptotic object has gsm}) we know
there exists $\rho>0$ such that 
$$
\|v\|_{L^2_{\{|y|<\rho\}}} > \|Q\|_{L^2} - \ep.
$$
Since, by Corollary \ref{modifiedcwlowercor},  the $\sigma_n$ go to
$\infty$ at least as fast as $(T^* - t_n)^{-s/2}$, eventually $\rho <
(T^* - t_n )^{s/2} \gamma(T^* - t_n) \sigma_n$, where $\gamma(z)$ is
a positive function satisfying $\gamma (z) \uparrow \infty$ as $z \downarrow 0$.
Thus,
\begin{equation*}
\begin{split}
\lim_{n\to\infty}\|u(t_n)\|_{L^2_{\{|x|<(T^* - t_n )^{s/2} \gamma(T^*
    - t_n)\}}} &=
\lim_{n\to\infty}\|v_n\|_{L^2_{\{|y|< (T^* - t_n )^{s/2} \gamma(T^*
    - t_n) \sigma_n\}} }\\
                                         &\ge\lim_{n\to\infty}\|v_n\|_{L^2_{
                                         \{|y|<\rho\}}}\\
                                         &\ge \|v\|_{L^2_{\{|y|<\rho\}}}\\
                                         &> \|Q\|_{L^2} - \ep.
\end{split}
\end{equation*}
This concludes the proof of Theorem \ref{main result} under the
assumptions that  Proposition
\ref{umbrella} and Corollary \ref{modifiedcwlowercor} hold true.
\end{proof}

\section{local-in-time theory}
\label{localtheorysection}
In this section, we adapt the arguments in \cite{CKSTT:MRL} to prove
an almost conservation property for $E[I_N u]$ which plays a central
role in the proof of Proposition \ref{umbrella}. We
begin by revisiting the Strichartz estimates and the classical 
proof \cite{CW} of local
well-posedness of \eqref{nls} for $u_0 \in H^s, ~ s>0$. 
We then explain a modification 
of the $H^s$ local well-posedness result in which the $H^s$-norm of 
$u_0$  is 
replaced by $\| I_N \langle D \rangle u_0 \|_{L^2}$. The modified local
well-posedness result provides 
the space-time control used to prove the almost conservation
of the modified energy $E[I_N u]$.

\subsection{Strichartz estimates}

We recall the classical Strichartz estimates  \cite{GV} 
for the Schr\"odinger group
$e^{it\Delta}$ on $\R_t \times \R_x^2$ (see also \cite{KT} for a
unified presentation).

The ordered exponent pairs $(q,r)$ are
{\it{admissible}} if $\frac{2}{q} + \frac{2}{r} = 1, ~ 2 < q$. Note that
$(3,6)$ and $(\infty, 2)$ are admissible. 
We define the {\it{Strichartz norm}} of functions $u:
[0,T] \times \R^2 \rightarrow \C$,
\begin{equation}
  \label{StrichNorm}
  \| u \|_{S^0_T} = \sup_{(q,r) ~{\mbox{admissible}}} \| u \|_{L^q_{t
  \in [0,T]} L^r_{x \in \R^2}}.
\end{equation}
We will use the shorthand notation $L^q_T$ to denote $L^q_{t
  \in [0,T]}$ and $L^p_x$ for $L^p (\R^2_x)$.
The H\"older dual exponent of $q$ is denoted $q'$, so $\frac{1}{q}
+ \frac{1}{q'} =1 $. The Strichartz estimates may be expressed as
follows:
\begin{equation}
  \label{StrichEst}
  \| u \|_{S^0_T} \leq \| u(0) \|_{L^2} + \| ( i \partial_t + \Delta)
  u \|_{L^{q'}_{T} L^{r'}_{x \in \R^2}}
\end{equation}
where $(q,r)$ is \emph{any} admissible exponent pair.

The smoothing properties underlying our proof of the almost
conservation of $E[I_N u]$ given below requires a careful control of
the interaction between high and low frequency parts of the
solution. Linear Strichartz estimates are not sufficient for this
purpose. They are complimented by bilinear Strichartz estimates that
were introduced in \cite{BRefine} 
and revisited in \cite{CKSTT:gopher}.
Let $D^\alpha$ denote the operator defined by 
 $\widehat{D^\alpha u} (\xi ) = |\xi|^\alpha \widehat{u} 
(\xi).$ Similarly, $\langle D \rangle^\alpha$ denotes the operator
defined via $\widehat{\langle D \rangle^\alpha u} (\xi ) = (1+ |\xi|^2)^{\frac{\alpha}{2}} \widehat{u}(\xi).$
We recall the following {\it{a priori bilinear Strichartz estimate}}
\cite{BRefine} as expressed in Lemma 3.4 of 
\cite{CKSTT:gopher}: For all $\delta
> 0$ and any $u,v: [0,T] \times \R^2 \longmapsto \C$,
\begin{eqnarray}
  \label{BRS} \ \ \ \ \
  \| u v \|_{L^2_T L^2_x} \leq C(\delta) & \left ( \| D^{-\half + \delta} u(0)
  \|_{L^2_x} + \| D^{-\half + \delta} (i \partial_t + \Delta) u
  \|_{L^1_T L^2_x}\right )  \\ & \times \left ( \| D^{\half - \delta} v(0)
  \|_{L^2_x} + \| D^{\half - \delta} (i \partial_t + \Delta) v
  \|_{L^1_T L^2_x}\right ). \notag
\end{eqnarray}
Note that, in the inequality \eqref{BRS},  a different amount of
regularity is required for $u$ and $v$. In the following analysis, we
will use it with $u,v$ being the projection of an $H^s$ solution of
\eqref{nls} onto different frequency regimes.

\subsection{Standard $H^s$ Local well-posedness}

We revisit the proof of local well-posedness of \eqref{nls} for
initial data in $H^s, ~ s>0$, extracting the features
of the local-in-time theory required in the proof of the almost
conservation property underlying our proof of Theorem 
\ref{main result} and Proposition \ref{umbrella}. 

\begin{proposition}[$H^s$-LWP \cite{CW}] 
\label{lwpprop}For $u_0 \in H^s (\R^2), ~ s>0$, the
  evolution $u_0 \longmapsto u(t)$ is well-posed on the time interval
  $[0, T_{lwp}]$ with
  \begin{equation}
    \label{lifetime}
    T_{lwp} = c_0 \| \langle D \rangle^s u_0 \|_{L^2}^{-\frac{2}{s}},
  \end{equation}
for a constant $c_0$ and
  \begin{equation}
    \label{doublingcontrol}
    \| \langle D \rangle^s u \|_{S^0_{T_{lwp}}} \leq 2 \| \langle D
    \rangle^s u_0 \|_{L^2_x}.
  \end{equation}
\end{proposition}

\begin{proof}
  The initial value problem \eqref{nls} is
  equivalent, by Duhamel's formula, to solving the integral equation
  \begin{equation}
    \label{Duhamel}
    u(t) = e^{it \Delta} u_0 + i \int_0^t e^{i (t - t') \Delta} (
    |u|^2 u (t') ) dt'.
  \end{equation}
For a given function $u$, denote the right-hand  side of \eqref{Duhamel}
by $\Phi_{u_0} [u]$. We prove that $\Phi_{u_0} [ \cdot ]$ is a
contraction mapping on the ball
\begin{equation}
  \label{Ball}
  B_{L^3_T L^6_x} ( \rho ) = \{ u: [0,T] \times \R^2 \longmapsto \C
 ~\big{|}~ \| u \|_{L^3_T L^6_x } < \rho \},  
\end{equation}
for sufficiently small $\rho$. In fact, we will show that
$\|\Phi_{u_0} [u ] - \Phi_{u_0}[v] \|_{S^0_T} < \half \| u - v
\|_{L^3_T L^6_x}$ provided $T$ is chosen small enough in terms of the
$H^s$ size of the initial data.
We write $\Phi_{u_0} [u ] - \Phi_{u_0}
[v]$, observe the cancellation of the linear pieces, apply the $S^0_T$
norm, and use \eqref{StrichEst} to obtain
\begin{equation}
  \label{poststrich}
  \|\Phi_{u_0} [u ] - \Phi_{u_0} [v] \|_{S^0_T} \leq C \| |u|^2 u -
  |v|^2 v \|_{L^1_T L^2_x}.
\end{equation}
(We have made the choice $q'=1$ and $r' = 2$ in \eqref{StrichEst} to
simplify the analysis below which emerges from the appearance of this
norm in \eqref{BRS}. This choice makes it convenient to perform the
fixed point argument in $L^3_T L^6_x.$)

By H\"older's inequality and simple algebra,
\begin{equation*}
   \| |u|^2 u -
  |v|^2 v \|_{L^1_T L^2_x} \leq C (\| u \|_{L^3_T L^6_x}^2 + \| v
  \|_{L^3_T L^6_x }^2 ) \| u - v \|_{L^3_T L^6_x} \leq C \rho^2 \|u -
  v\|_{L^3_T L^6_x }.
\end{equation*}
This proves $\Phi_{u_0} [\cdot ]$ is indeed a contraction mapping on $B
(\rho)$ if $\rho < \rho_0$, where $\rho_0$ is an explicit constant.

Thus, the sequence of Duhamel iterates $\{ u^j \}$, defined by the
recursion
\begin{equation}
  \label{iterates}
\left\{
\begin{matrix}
  u^0 (t, x) = e^{i t \Delta } u_0 (x) \\
  u^{j+1} (t) = e^{it \Delta} u_0 + i \int_0^t e^{i (t - t') \Delta} (
  |u^{j} |^2 u^j (t') ) dt',
\end{matrix}
\right.
\end{equation}
converges geometrically to the solution $u$ provided that 
\begin{equation}
  \label{initalsize}
  \| e^{it \Delta } u_0 \|_{L^3_T L^6_x }  < \rho_0.
\end{equation}
By H\"older's inequality in time and Sobolev's inequality in space, 
\begin{equation}
  \label{HolderSob}
  \| u \|_{L^3_T L^6_x }  \leq C T^{\frac{1}{3} - \frac{1}{q}} \|
  \langle D \rangle^s u \|_{L^q_T L^r_x}
\end{equation}
with $q > 3$ and $\frac{2}{r} = \frac{2}{6} + s $. We choose $q$ so
that $(q,r)$ is admissible, {\it i.e.}
$\frac{2}{q} + \frac{2}{r} = 1$.  This gives $\frac{1}{3} -
  \frac{1}{q} = \frac{s}{2}.$
Applying (\ref{HolderSob}) and then \eqref{StrichEst}, we have
\begin{equation}
\begin{split}
\label{timecontrol}
  \| e^{i t \Delta} u_0 \|_{L^3_T L^6_x}  \leq & T^{\frac{s}{2}} \|  e^{i t \Delta}\langle D \rangle^s u_0 \|_{L^q_T L^r_x}\\
                                          \leq & T^{\frac{s}{2}} \| \langle D \rangle^s u_0 \|_{L^2}.
\end{split}
\end{equation}

We choose $T = T_{lwp}$ so that the right-hand side of \eqref{timecontrol}
is less than $\rho_0$, namely so that \eqref{lifetime} holds. 
The Strichartz estimate implies that the zeroth iterate satisfies $\|
u^0 \|_{S^0_{T_{lwp}}} \leq C \|u_0 \|_{H^s}$. By the geometric
convergence of the iterates, we have that
\begin{equation}
\label{spacetimecontrol}
  \| u \|_{S^0_{T_{lwp}}}  \leq C\| u_0 \|_{H^s}
\end{equation}
and
\begin{equation}
  \label{36localsmall}
  \| u \|_{L^3_{T_{lwp}} L^6_x} \leq C \rho_0.
\end{equation}

A posteriori, we can revisit the estimate for the solution of
\eqref{Duhamel} and use the Leibnitz rule
for fractional derivatives to prove the persistence of regularity
property 
\begin{equation}
  \label{persistence}
  \| \langle D \rangle^s u \|_{S^0_T }  \leq 2 \| \langle
  D \rangle^s u_0 \|_{L^2}.
\end{equation}
This follows from the Strichartz estimate on the linear term in
\eqref{Duhamel} and the estimate
\begin{equation}
  \label{HolderLeibnitz}
  \| \langle D \rangle^s (|u|^2 u )  \|_{L^1_T L^2_x} \leq C \|u
  \|_{L^3_T L^6_x}^2  \| \langle D \rangle^s u \|_{L^3_T L^6_x } \leq
 C \rho_0^2  \| \langle D \rangle^s u \|_{S^0_T}
\end{equation}
where, in the last step, we used \eqref{36localsmall}. The first
inequality in \eqref{HolderLeibnitz} follows, for example,  from Proposition 1.1 on
page 105 of \cite{TaylorTools}.

\end{proof}

\begin{corollary}[\cite{CW}]
  If $H^s \ni u_0 \longmapsto u(t)$ with $s>0$ solves \eqref{nls} for all $t$
  near enough to $T^*$ in
  the maximal finite interval of existence $[0, T^*)$ then
  \begin{equation}
    \label{cwlowerbound}
C    (T^* - t )^{-\frac{s}{2}} \leq \| \langle D \rangle^s u(t) \|_{L^2_x}.
  \end{equation}
\end{corollary}
  Note that the estimate \eqref{HolderLeibnitz} implies that, when we
  restrict attention to 
  solutions of \eqref{nls} with $t \in [0, T_{lwp}]$, 
we can 
  essentially ignore the contributions involving the 
  inhomogeneous terms appearing on the right side of \eqref{BRS}
  since these terms are controlled by the corresponding homogeneous
  terms. 
As a consequence of this and the bilinear Strichartz estimate,  
we prove the following lemma. We define $\proj_{N_j}$
to be the Littlewood-Paley projection operator to frequencies of size
  $N_j \in 2^{\mathbb{N}}$, i.e.
$\widehat{\proj_{N_j} f } ( \xi ) = \chi_{\{\frac12 N_j < \xi < 2
  N_j\}} \hat{f} (\xi)$.  

\begin{lemma}
  \label{bilinearstrich}
Suppose $u$ solves \eqref{nls} on the time interval $[0,
T_{lwp}]$. Let $u_{j} = \proj_{N_j} u$, for $j=1,2$ with $N_1 >
N_2$. Then
\begin{equation}
  \label{2linearbrs}
  \| u_{1} u_{2} \|_{L^2_{T_{lwp}} L^2_x} \leq C
  \left(\frac{N_2}{N_1} \right)^{\half -} \| u \|_{S^0_{T_{lwp}}}^2,
\end{equation}
The estimate \eqref{2linearbrs} is
also valid if $u_j$ is replaced by $\overline{u_j}$.
\end{lemma}

This bilinear smoothing property of solutions of \eqref{nls} is the
key device underpinning the proof of the almost conservation of the
modified energy in Proposition \ref{ACProp}.

\subsection{Modified $H^s$ Local Well-posedness}

\begin{proposition}[Modified $H^s$-LWP]
\label{modifiedHsLWP} For $u_0 \in H^s,~s>0$, the  \eqref{nls}
evolution $u_0 \longmapsto u(t)$ is well-posed on the time interval
  $[0, \widetilde{T}_{lwp}]$ with
  \begin{equation}
    \label{modifiedlifetime}
    \widetilde{T}_{lwp} = c_0 \| I_N \nabla u_0 \|_{L^2_x}^{-\frac{2}{s}},
  \end{equation}
  \begin{equation}
    \label{modifieddoublingcontrol}
    \| I_N \langle D \rangle u \|_{S^0_{\widetilde{T}_{lwp}}} \leq 2
    \| I_N \langle D \rangle    u_0 \|_{L^2_x}.
  \end{equation}
\end{proposition}

\begin{proof}
The proof is a modification of the arguments used to prove Proposition 
\ref{lwpprop}. Note first that the estimates \eqref{timecontrol} and
\eqref{easy I properties} may be combined to give 
\begin{equation}
  \label{modifedtimecontrol}
  {\| e^{it \Delta} u_0 \|}_{L^3_T L^6_x} \leq C T^{\frac{s}{2}} {{\|
  I_N \langle D \rangle u_0 \|}_{L^2}}.
\end{equation}
Thus, the previous analysis produces a solution $u$ of \eqref{nls}
satisfying \eqref{36localsmall} and 
\begin{equation}
  \label{modifiedspacetimecontrol}
  \| u \|_{S^0_{\widetilde{T}_{lwp}}} \leq C \| I_N \langle D \rangle
  u_0 \|_{L^2}
\end{equation}
provided that we choose ${\widetilde{T}}_{lwp}$ as in \eqref{modifiedlifetime}.

Next, we turn our attention toward the space-time regularity estimate
\eqref{modifieddoublingcontrol}. Since $e^{it\Delta}$ does not affect
the magnitude of Fourier coefficients, \eqref{modifieddoublingcontrol}
is clearly valid for the  linear term in \eqref{Duhamel}. Since
the spaces appearing on both sides of the trilinear estimate
\eqref{HolderLeibnitz} are translation invariant and (the first
inequality in)
\eqref{HolderLeibnitz} is valid, the modified estimate
\eqref{modifieddoublingcontrol} follows directly from the
interpolation lemma (Lemma 12.1 on page 108) in \cite{CKSTT:gKdV}.
\end{proof}

\begin{remark}
A modification of \eqref{2linearbrs} follows using the spacetime
control \eqref{modifieddoublingcontrol}: For $N_1 \geq N_2$ and for
solutions $u$ of \eqref{nls},
\begin{equation}
\label{modifiedBRS}
\| I \langle D \rangle u_{N_1} I \langle D \rangle u_{N_2}
\|_{L^2_{\widetilde{T}} L^2_x} \leq C \left( \frac{N_2}{N_1} \right)^{\half-}
 \| I \langle D \rangle u \|^2_{S^0_{\widetilde{T}}}.
\end{equation}
\end{remark}

\begin{corollary}
\label{modifiedcwlowercor}
  If $H^s \ni u_0 \longmapsto u(t)$ with $s>0$ solves \eqref{nls} for all $t$
  near enough to $T^*$ in
  the maximal finite time interval of existence $[0,T^*)$
  \begin{equation}
    \label{modifiedcwlowerbound}
    C (T^* - t )^{-\frac{s}{2}} \leq \| I_N \langle D \rangle u(t) \|_{L^2_x}.
  \end{equation}
\end{corollary}

Since we are studying here finite time blowup solutions, we will
sometimes implicitly assume that $\| I_N \nabla u(t) \|_{L^2} >
1$. 

\subsection{Almost conservation law for the modified energy}

\begin{proposition}
  \label{ACProp}
If $H^s \ni u_0 \longmapsto u(t)$ with $s>0$ solves \eqref{nls} for all $t \in [0,
\widetilde{T}_{lwp}]$ then
\begin{equation}
  \label{ACest}
  \sup_{t \in [0, \widetilde{T}_{lwp}]} |E[ I_N u(t) ]| \leq
 |E[I_N u(0) ]|
  + C N^{-\alpha_4 } \| I_N \langle D \rangle u(0) \|_{L^2_x}^4 + C N^{-\alpha_6}
  \|I_N \langle D \rangle u(0) \|_{L^2_x}^6,
\end{equation}
with $\alpha_4 = \frac{3}{2}-$ and $\alpha_6 = 2-$.
\end{proposition}

\begin{proof}
  We adapt arguments from \cite{CKSTT:MRL} in which a similar result
  is proved using local well-posedness theory in the weighted 
$X_{s,b}$ spaces. We recall that the parameter $N$ 
refers to the operator $I_N$ defined in
(\ref{mdefined}). In light of
  \eqref{modifieddoublingcontrol}, it
  suffices to control the {\it{energy increment}} $| E[I_N u(t) ] -
  E[I_N u(0) ] |$ for $ t \in [0, \widetilde{T}_{lwp}]$ in terms of
  $\| I_N \langle D \rangle u \|_{S^0_{\widetilde{T}_{lwp}}}$. We define the set
$*_n = \{ (\xi_1, \ldots, \xi_n) : \Sigma \xi_i = 0 \}$. 
 An explicit calculation  of $\partial_t E[I_N u]$ 
(carried out in detail in Section 3 of  \cite{CKSTT:MRL})
 reveals that $| E[I_N u(t) ] - E[I_N u(0)]|$ is
  controlled by the sum of the two space-time integrals:

  \begin{equation}
    \label{4linearexp}
E_1= \left|    \int_0^t \int_{*_4}  [ 1 - \frac{m(\xi_1 )}{m(\xi_2 ) m(\xi_3)
    m(\xi_4)} ] \widehat{ \Delta I \overline{u}} (\xi_1 ) \widehat{I u
    } (\xi_2 ) \widehat{ I \overline{u}} (\xi_3 ) \widehat{I u } (\xi_4)\right|,
  \end{equation}
and
  \begin{equation}
    \label{6linearexp}
E_2 =\left|     \int_0^t \int\limits_{\small{\tilde{\xi_1} + \xi_4 + \xi_5 + \xi_6 = 0}}
  [ 1 - \frac{m({\tilde{\xi_1}} )}{m(\xi_4 ) m(\xi_5)
    m(\xi_6)} ] {\widehat{ I ( |u|^2 \overline{u} )}} ( \tilde{\xi_1} )  \widehat{I u
    } (\xi_4 ) \widehat{ I \overline{u}} (\xi_5 ) \widehat{I u } (\xi_6)\right|.
  \end{equation}

We estimate the 4-linear expression \eqref{4linearexp} first. Let
$u_{N_j}$ denote $\proj_{N_j} u$. When $\xi_j$ is dyadically localized
to $\{ |\xi | \thicksim N_j \}$ we will write $m_j$ to denote $m(\xi_j)$.

\begin{lemma}
  \label{4linearlemma} 
If $H^s \ni u_0
  \longmapsto u(t)$ with $s>0$ solves
  \eqref{nls} for all $t \in [0 , \widetilde{T}_{lwp} ]$ then
  \begin{equation}
   \label{4linearest}
\begin{split}   
&\left| \int_0^t \int_{*_4}  [ 1 - \frac{m_1}{m_2 m_3 m_4 }]
\left ( \widehat{ \Delta I \overline{u}_{N_1}} (\xi_1 ) \widehat{I
u_{N_2}
    } (\xi_2 ) \widehat{ I \overline{u}_{N_3}} (\xi_3 ) \widehat{I
      u_{N_4} } (\xi_4)\right ) \right| \\
& \qquad \leq N^{-\alpha_4} \| I
    \langle D \rangle u \|^4_{S^0_{\widetilde{T}_{lwp}}} \prod_{j=1}^4 N_j^{0-},
\end{split}
  \end{equation}
with $\alpha_4 = \frac{3}{2} -$,
\end{lemma}
\begin{proof}

The analysis which follows will not rely upon the complex conjugate
structure in the left-side of \eqref{4linearest}. Thus, there is
symmetry under the interchange of the indices 2,3,4. We may therefore
assume that $N_2 \geq N_3 \geq N_4.$

{\bf{Case 1.}} $N \gg N_2$. On the convolution hypersurface $*_4$,
this forces $N_1 \ll N$ as well, so the multiplier $ [ 1 -
\frac{m_1}{m_2 m_3 m_4 }] = 0$ and the expression to be bounded vanishes.

{\bf{Case 2.}} $N_2 \gtrsim N \gg N_3 \geq N_4$. This forces $N_1
\thicksim N_2$ on $*_4$. By the mean value theorem and simple algebra
\begin{equation}
  \label{MVTmultest}
  \left|  [ 1 - \frac{m_1}{m_2 m_3 m_4 }] \right| = \left| [ \frac{m(\xi_2) - m(\xi_2
     + \xi_3 + \xi_4)}{m(\xi_2) } ]\right| \leq \left| \frac{\nabla m(\xi_2 )
     \cdot (\xi_3 + \xi_4 )}{m(\xi_2)} \right| \lesssim \frac{N_3}{N_2}.
\end{equation}

We use the frequency localization of the $u_{N_j}$ to renormalize the
derivatives and multipliers to arrange for the appearance of $I
\langle D \rangle u_{N_j}$. In the case under consideration, using
\eqref{MVTmultest} and the frequency localizations, the left-side of 
\eqref{4linearest} is controlled
by 
\begin{equation*}
  \frac{N_3}{N_2} N_1 ( N_2 \langle  N_3 \rangle \langle N_4 \rangle)^{-1}  \left| \int_0^t \int_{*_4}
  \prod_{j=1}^4  I \langle D \rangle u_{N_j} (\xi_j ) \right|. 
\end{equation*}
We apply Cauchy-Schwarz to obtain the bound
\begin{equation*}
  \frac{N_3}{N_2} N_1 ( N_2 \langle N_3 \rangle  \langle N_4 \rangle  )^{-1} \| I \langle D \rangle
  u_{N_1}  I \langle D \rangle u_{N_3} \|_{L^2_{{\widetilde{T}_{lwp}} L^2_x}} \| I \langle D \rangle
  u_{N_2}  I \langle D \rangle u_{N_4} \|_{L^2_{\widetilde{T}_{lwp}} L^2_x}.
\end{equation*}
By \eqref{modifiedBRS}, we control by
\begin{equation*}
  \frac{N_3}{N_2}  N_1 (N_2 \langle N_3 \rangle  \langle N_4 \rangle )^{-1}  \left( \frac{N_3 N_4 }{N_1
  N_2} \right)^{\half -} \| I \langle D \rangle u \|^4_{S^0_{\widetilde{T}_{lwp}}}
\end{equation*}
which simplifies  to give the bound
\begin{equation}
  \label{Case2final}
  N^{-\frac{3}{2} +} (N_1 N_2 \langle N_3 \rangle \langle N_4 \rangle 
)^{0-} \| I \langle D \rangle
  u \|^4_{S^0_{\widetilde{T}_{lwp}}}.
\end{equation}

{\bf{Case 3.}} $N_2 \geq N_3 \gtrsim N.$ In this case, we use the
trivial multiplier bound
\begin{equation}
  \label{trivialmultest}
  \left|  [ 1 - \frac{m_1}{m_2 m_3 m_4 }] \right| \leq  \frac{m_1}{m_2 m_3 m_4 }.
\end{equation}
Again, we pull the multiplier out and estimate the remaining integral
using Lemma \ref{bilinearstrich}.

{\bf{Case 3a.}} $N_1 \thicksim N_2 \geq N_3 \gtrsim N$. We bound 
$E_1$ here by renormalizing the derivatives and multiplier, then pairing $u_{N_1} u_{N_3}$ and
$u_{N_2} u_{N_4}$ and using Lemma \eqref{bilinearstrich} again:
\begin{equation*}
   \frac{m_1}{m_2 m_3 m_4 } \left( \frac{N_3 N_4 }{N_1 N_2}
   \right)^{\half - } N_1 (N_2 N_3 \langle N_4 \rangle )^{-1}  \| I
  \langle D \rangle u \|^4_{S^0_{\widetilde{T}_{lwp}}}.
\end{equation*}
We reexpress this bound as
\begin{equation*}
  \frac{m_1}{m_2 N_2^\half m_3 N_3^\half m_4 \langle N_4 \rangle^\half
  N_1^\half}  \| I \langle D \rangle
  u \|^4_{S^0_{\widetilde{T}_{lwp}}}.
\end{equation*}
Since $m(x)$ is bounded from above by 1 and $m(x) \langle x \rangle ^\half$ is 
nondecreasing and bounded from below by 1, 
this is bounded by
\begin{equation}
  \label{Case3afinal}
  \frac{1}{N N_1^\half}  \| I \langle D \rangle
  u \|^4_{S^0_{\widetilde{T}_{lwp}}} \leq N^{-\frac{3}{2} + } (N_1
  N_2 N_3 \langle N_4 \rangle )^{0-}  \| I \langle D \rangle
  u \|^4_{S^0_{\widetilde{T}_{lwp}}}.
\end{equation}
{\bf{Case 3b.}} $N_2 \thicksim N_3 \geq N$. A similar analysis leads
to the bound
\begin{align*}
  &  \frac{m_1}{m_2 m_3 m_4 } \left( \frac{N_1 N_4 }{N_3 N_2}
   \right)^\half  N_1 (N_2
   N_3 \langle N_4 \rangle)^{-1}   \| I \langle D \rangle
  u \|^4_{S^0_{\widetilde{T}_{lwp}}} \\
& \leq \frac{1}{N^\half N_2}  \| I \langle D \rangle
  u_{N_j} \|^4_{S^0_{\widetilde{T}_{lwp}}} \\
& \leq N^{-\frac{3}{2} + } (N_1 N_2 N_3 \langle N_4 \rangle )^{0-} \| I \langle D \rangle
  u_{N_j} \|^4_{S^0_{\widetilde{T}_{lwp}}}.
\end{align*}

\end{proof}

A related case-by-case analysis combined with a trilinear estimate
establishes the required 6-linear estimate for \eqref{6linearexp}. We
write $m_{123}$ to denote $m(\xi_1 + \xi_2 + \xi+_3)$ and use
$N_{123}$ to denote the (dyadic)
size of $\xi_1 + \xi_2 + \xi_3$. We will also use the similarly
defined notation $m_{456}$. Note that $N_{123}$ could be much
smaller than $N_1, N_2,$ or $N_3$.

\begin{lemma}
  \label{6linearlemma}  If $H^s \ni u_0 \longmapsto u(t)$ with $s>0$ solves
  \eqref{nls} for all $ ~t \in [0 , \widetilde{T}_{lwp} ]$ then
  \begin{equation*}
    \int_0^t \int_{*_6}  \left[ 1 - \frac{m_{123}}{m_4 m_5 m_6} \right]
         m_{123} \left( {\widehat{\overline{u}}}_{N_1} ( \xi_1 )
    {\widehat{u}}_{N_2} (\xi_2 ) {\widehat{ \overline{u}}}_{N_3} (
    \xi_3 ) \right)   {\widehat{Iu}}_{N_4} (\xi_4)
 {\widehat{{ I \overline{u}}}}_{N_5} (\xi_5 )
 { \widehat{I u}}_{N_6} (\xi_6 )
\end{equation*}
\begin{equation}
 \label{6linearest}
 \leq N^{-\alpha_6} \| I \langle D \rangle u
\|^6_{S^0_{\widetilde{T}_{lwp}}} \prod_{j=1}^6 N_j^{0-}  
  \end{equation}
with $\alpha_6 = 2-$.
\end{lemma}

\begin{proof}
  We carry out a case-by-case analysis. By symmetry (since we will not
  use the complex conjugate structure), we may assume $N_4 \geq N_5
  \geq N_6$.

{\bf{Case 1.}} $N \gg N_4$. On $*_6$, this forces $N_{123} \thicksim
N_4$ so the multiplier $ \left[ 1 - \frac{m_{123}}{m_4 m_5 m_6} \right]$ 
  vanishes.

{\bf{Case 2.}} $N_4 \gtrsim N \geq N_5$. On $*_6$, $N_{123} \thicksim
N_4$ in this case. By the mean value theorem,
\begin{equation*}
  \left|  \left[ 1 - \frac{m_{123}}{m_4 m_5 m_6} \right] \right| = \left| \frac{ m_4 - m_{456}
  }{m_4} \right| \leq \frac{N_5}{N_4}.
\end{equation*}

Applying this multiplier bound and the Cauchy-Schwarz inequality 
 to the integral in
\eqref{6linearest} gives the bound
\begin{equation*}
(N_4 \langle N_5 \rangle )^{-1} \frac{N_5}{N_4}  \| \proj_{N_{123}} I ( {\overline{u}}_{N_1} u_{N_2}
{\overline{u}}_{N_3} ) I u_{N_6} \|_{L^2_{Tx}} \| I \langle D \rangle
u_{N_4}  I \langle D \rangle u_{N_5}
\|_{L^2_{Tx}} .
\end{equation*}
By H\"older's inequality  and Lemma \ref{bilinearstrich}, we control the
above expression  by
\begin{equation}
  \label{saveme}
  (N_4 \langle N_5 \rangle )^{-1} 
 \frac{N_5}{N_4}    \| \proj_{N_{123}} I ( {\overline{u}}_{N_1} u_{N_2}
{\overline{u}}_{N_3} )\|_{L^2_{Tx}} \|  I u_{N_6} \|_{L^\infty_{Tx}}
 \left( \frac{N_5}{N_4}
\right)^\half  \| I  \langle D \rangle u  \|_{S^0_{{\widetilde{T}}_{lwp}}}  \| I \langle D \rangle
u    \|_{S^0_{{\widetilde{T}}_{lwp}}} .
\end{equation}
By Sobolev's inequality on
functions with frequency support localized to a dyadic shell on
$R^2_x$,
\begin{equation}
  \label{2dstep}
  \| I u_{N_6} \|_{L^\infty_{Tx} } \leq \| I \langle D \rangle u_{N_6}
  \|_{L^\infty_T L^2_x} .
\end{equation}
In order to continue the proof of Lemma \ref{6linearlemma}, we need
an estimate of the term $\| \proj_{N_{123}} I ( {\overline{u}}_{N_1} u_{N_2}
{\overline{u}}_{N_3} )\|_{L^2_{Tx}}$. 
This is the purpose of the next lemma.
Let $N_1^\dagger  \geq N_2^\dagger  \geq N_3^\dagger  $ denote the
decreasing rearrangement of $N_1, N_2 , N_3$. 

\begin{lemma}
 \label{trilinearlemma}
 \begin{equation}
   \label{trilinearest}
   \| \proj_{N_{123}} I ( {\overline{u}}_{N_1} u_{N_2}
{\overline{u}}_{N_3} )\|_{L^2_{Tx}}  \leq 
\langle N_1^\dagger \rangle^{-\frac{1}{2}}
\prod_{j=1}^3  \| I \langle D \rangle
  u_{N_j} \|_{S^0_{\widetilde{T}_{lwp}}}.
 \end{equation}
\end{lemma}

\begin{proof} 
Again, we will not use the complex conjugate structure so we may
assume that $N_1 \geq N_2 \geq N_3$. The projection $\proj_{N_{123}}$
allows us to control the left-hand side of \eqref{trilinearest} by
\begin{equation*}
  m_{123} \| \overline{u}_{N_1}  u_{N_2}  \overline{u}_{N_3}
  \|_{L^2_{Tx}}  \leq m_{123}  \| u_{N_1} \|_{L^4_{Tx}}  \| u_{N_2}
  \|_{L^4_{Tx}} \|u_{N_3} \|_{L^\infty_{Tx}}
\end{equation*}
where we used H\"older's inequality.

We renormalize the
derivatives and use a (dyadically localized) Sobolev inequality
 as in \eqref{2dstep} to get
the bound
\begin{equation*}
  \frac{m_{123}}{m_1 m_2 m_3 }  \langle N_1 \rangle^{-1}  \langle N_2
  \rangle^{-1}  \| I \langle D \rangle u_{N_1} \|_{L^4_{Tx}}
\| I \langle D \rangle u_{N_2} \|_{L^4_{Tx}}    \| I \langle D \rangle u_{N_3} \|_{L^\infty_T
L^2_x} .
\end{equation*}
Since the norms appearing in the above expression are
admissible we focus our attention upon the prefactor and bound with
the expression
\begin{equation*}
  \frac{ m_{123}  N_3^\half }{\langle N_1 \rangle^\half \langle N_2
 \rangle^\half m_1 \langle N_1 \rangle^\half m_2
 \langle N_2 \rangle^\half m_3 \langle N_3 \rangle^\half  }
 \prod_{j=1}^3   \| I \langle D \rangle
  u_{N_j} \|_{S^0_{\widetilde{T}_{lwp}}}.
\end{equation*}
Since $m(x) \leq 1, ~ m(x) \langle x \rangle^{\half}$ is
nondecreasing, and $N_3 \leq N_2$, this proves \eqref{trilinearest}.
\end{proof} 

We use \eqref{trilinearest} and \eqref{2dstep} on \eqref{saveme} to
complete the Case 2 analysis.  The left-hand side of
 \eqref{6linearest} is bounded by 
\begin{eqnarray}
  \label{Case2final6}
& (N_4 \langle N_5 \rangle )^{-1} \left( \frac{N_5}{N_4} \right)^{\frac{3}{2}}
  { \langle N_{1}^\dagger \rangle}^{-\half}    \| I \langle D \rangle
  u \|^6_{S^0_{\widetilde{T}_{lwp}}}\\
& \leq \frac{N_5^{\frac{1}{2}}}{N_4^{\frac{5}{2}}}
(N_1^\dagger)^{-\half}  \| I \langle D \rangle
  u \|^6_{S^0_{\widetilde{T}_{lwp}}}  \\
& \leq N^{-2+} (N_1 \dots \langle N_5 \rangle \langle N_6 \rangle )^{0 - } \prod_{j=1}^6   \| I \langle D \rangle
  u \|^6_{S^0_{\widetilde{T}_{lwp}}}.
\end{eqnarray}

{\bf{Case 3.}}  $N_4 \geq N_5 \geq N.$  We will use the trivial
multiplier estimate
\begin{equation*}
     \left |  \left[ 1 - \frac{m_{123}}{m_4 m_5 m_6} \right] \right| \leq  \frac{m_{123}}{m_4 m_5 m_6}.
\end{equation*}
Familiar steps lead to the bound
\begin{equation*}
  \frac{m_{123}}{m_4 m_5 m_6}  \langle N_1^\dagger \rangle^{-\half}  \left(
  \frac{N_5}{N_4} \right)^{\frac{1}{2}} (N_4 N_5)^{-1} 
\| I \langle D \rangle   u \|^6_{S^0_{\widetilde{T}_{lwp}}}
\end{equation*}
which we recast as
\begin{equation}
  \label{lookatme}
  \frac{m_{123}}{m_4 N_4^{\half} m_5 N_5^{\half} m_6 \langle
  N^\dagger_1 \rangle^{\half}
  N_4 }  \| I \langle D \rangle   u
  \|^6_{S^0_{\widetilde{T}_{lwp}}} .
\end{equation}

{\bf{Case 3a.}}  $N_6 \geq N$. We express the prefactor in
\eqref{lookatme} as
\begin{equation*}
  \frac{m_{123}  N_6^{\half}  }{ m_4 N_4^{\half}  m_5 N_5^{\half}  m_6
  N_6^{\half}  \langle N_1^{\dagger} \rangle^{\half} N_4 }\leq
  \frac{1}{N^{\frac{3}{2}} N_4^{\half}  }  \leq N^{-2+}  (N_1 \dots
  N_6 )^{0-}.
\end{equation*}

{\bf{Case 3b.}}  $N_6 \leq N$. Here we have $m_6 = 1$ so we can bound
the prefactor in \eqref{lookatme} by
\begin{equation*}
  \frac{1}{N N_4}  \leq N^{-2+}  (N_1 N_2 N_3 N_4)^{0-} .
\end{equation*}

\end{proof} 
 Since \eqref{4linearexp} and \eqref{6linearexp} control
 the increment in the modified energy, it suffices to prove
 appropriate bounds on these integrals to obtain \eqref{ACest}. Lemmas
 \eqref{4linearlemma}, \eqref{6linearlemma} provide estimates for
 dyadically localized contributions to the integrals
 \eqref{4linearexp}, \eqref{6linearexp}, respectively. Since the
 estimates \eqref{4linearest} and \eqref{6linearest} have the helpful
 decay factors $N_j^{0-}$, we can sum up the dyadic contributions, and
 apply Cauchy-Schwarz to complete the proof of \eqref{ACest}.

\end{proof} 

\section{Modified kinetic energy dominates modified total energy}
\label{umbrella lemma section}

In this section, we prove Proposition \ref{umbrella}.

\begin{proof}[Proof of Proposition \ref{umbrella}]
When $s \geq 1$, we set $I_{N(T)} = Identity$ by choosing $N(T) =
+ \infty$. Since the (unmodified) energy is conserved while the
kinetic energy blows up as time approaches $T^*$, the estimate
\eqref{energycontrol} is obvious with $p(s) = 0$. We therefore
restrict attention to $s \in (s_Q, 1)$, with $s_Q$ to be
determined at the end of the proof.

Fix $s \in (s_Q, 1)$ and $T$ near $T^*$. Let $N = N(T)$ (to be
chosen). Set $\delta = c_0 (\Sigma (T) )^{-\frac{2}{s}} > 0$ with
$c_0$ the small fixed constant in \eqref{modifiedlifetime}. Note that
$\delta$ is the time of local well-posedness guaranteed by Proposition
\ref{modifiedHsLWP} for initial data of size $\Sigma(T)$, which is the
largest value that the modified kinetic energy achieves up to time $T$.  Thus
the interval $[0,T]$ may be partitioned into $J = C \frac{T}{\delta} $
$\delta$-sized intervals on which the modified local well-posedness
result uniformly applies.  More precisely, 
$[0,T] = \bigcup_{j=1}^J I_j, ~ I_j = [t_j, t_{j+1}), ~ t_0=0, ~
t_{j+1} = t_j + c \delta$,
and we have at each $t_j$,
$$
\| I_N \langle D \rangle u(t_j)\|_{L^2} = \sigma(t_j ) \leq \Sigma(T).
$$
In addition, $\delta$ has been taken sufficiently small so
that we can apply the almost conservation law Proposition
\ref{ACProp} on each of the $I_j$.  
We now accumulate increments to the energy and have that
\begin{align*}
|  E[I_N u(T) ]| &\leq| E[I_N u(0)]| + C \frac{T^*}{\delta}  [
  N^{-\alpha_4} \Sigma^4 (T) + N^{-\alpha_6}  \Sigma^6 (T) ] \\
  & \leq N^{2(1-s)} \lambda (0) +  C \frac{T^*}{\delta} N^{-\alpha_4} \Sigma^4 (T)
  + C \frac{T^*}{\delta} N^{-\alpha_6} \Sigma^6 (T) .
\end{align*}
By the choice of $\delta$ we see, dismissing irrelevant
constants, that
\begin{equation*}
|  E[I_N u(T) ]| \lesssim N^{2(1-s)}  + N^{-\alpha_4} \Sigma^{4+
  \frac{2}{s}} (T)  + N^{-\alpha_6} \Sigma^{6 + \frac{2}{s} } (T).
\end{equation*}
Using \eqref{easy I properties}, we can switch from $\Sigma$ to $\Lambda$:
\begin{equation}
  \label{vladimir}
|  E[I_N u(T) ] |\lesssim N^{2(1-s)}  + N^{-\alpha_4 + (4 + \frac{2}{s}
  ) (1-s) }  \Lambda^{4 + \twoovers} (T) + N^{-\alpha_6 + (6 +
  \twoovers )(1-s) } \Lambda^{6 + \twoovers}  (T).
\end{equation}

We choose $N = N (\Lambda)$ so that the first and third terms
in \eqref{vladimir} give comparable contributions. 
A calculation reveals that the second term in \eqref{vladimir}
produces a smaller correction. 
Thus, with the choice
\begin{equation}
  \label{Nchoice}
  N = \Lambda^{\frac{6 + \twoovers}{\alpha_6 - (4 + \twoovers)(1-s)} },
\end{equation}
the Proposition \ref{umbrella} is established with
\begin{equation}
  \label{eq:pdefined}
  p(s) = \frac{6+\twoovers}{\alpha_6 - (4 + \twoovers )(1-s) } 2 (1-s).
\end{equation}
Note that $p(s) < 2$ reduces to 
 to a quadratic condition on $s$.  Specifically
$$
10s^2 + (\al_6 - 6)s - 4 > 0.
$$
For $\al_6 = 2-$ this yields
$$
s > s_Q = \frac{1}{5} + \frac{1}{5}\sqrt{11} \sim 0.863.
$$

\end{proof} 

\end{document}